\documentclass[final]{siamltex}

\usepackage{amssymb}
\usepackage{amsmath}
\usepackage{graphicx}
\usepackage{subfigure} 
\usepackage{color}

\newcommand \be {\begin{equation}}
\newcommand \ee {\end{equation}}
\newcommand \eps \epsilon  
\newcommand{\rf}[1]{(\ref{#1})}

\newcommand{\tA}{\widetilde{A}}

\newcommand{\bt}{b^{T}}
\newcommand{\tbt}{\widetilde{b}^{T}} 
\newcommand{\acapo}{\\[-0.3cm]} 

\begin{document} 
\title{High--Order Asymptotic--Preserving Methods
\\
for
\\
 fully nonlinear relaxation problems} 

\author{Sebastiano Boscarino\thanks{Department of
Mathematics and Computer Science, University of Catania, Via
A.Doria 6, 95125 Catania, Italy. ({\tt boscarino@dmi.unict.it, russo@dmi.unict.it})}\and
Philippe G. L{\scriptsize e}Floch\thanks{Laboratoire Jacques--Louis Lions \& Centre National de la Recherche Scientifique,
Universit\'e Pierre et Marie Curie (Paris 6), 4 Place Jussieu, 75252 Paris, France.
\newline 
 ({\tt contact@philippelefloch.org})}\and
Giovanni Russo$^\star$}

\maketitle

\begin{abstract} We study solutions to nonlinear hyperbolic systems with fully nonlinear relaxation terms in the limit of, both, infinitely stiff relaxation and arbitrary late time. In this limit, the dynamics is governed by effective systems of parabolic type, with possibly degenerate and/or fully nonlinear diffusion terms. For this class of problems, we develop here an implicit--explicit method based on Runge--Kutta discretization in time, and we use this method in order to investigate several examples of interest in compressible fluid dynamics. Importantly, we impose here a realistic stability condition on the time--step and we demonstrate that solutions in the hyperbolic--to--parabolic regime can be computed numerically with high robustness and accuracy, even in presence of fully nonlinear relaxation terms.  
\end{abstract}



\pagestyle{myheadings}
\thispagestyle{plain}
\markboth{S. BOSCARINO, P.G. L{\scriptsize E}FLOCH, AND G. RUSSO}{} 


\section{Introduction}

We consider nonlinear hyperbolic systems containing fully nonlinear relaxation terms when the relaxation is stiff relaxation and for 
arbitrary late time. The dynamics is asymptotically governed by effective systems which are parabolic type and may contain possibly
 degenerate and fully nonlinear diffusion terms. For such problems, we introduce implicit--explicit (IMEX) methods
which are based on Runge--Kutta (R--K) discretization in time, and we apply this method 
to investigate several compressible fluid models. Importantly, we impose here a realistic stability condition on the time--step and we demonstrate that solutions in the hyperbolic--to--parabolic regime can be computed numerically with high robustness and accuracy, even in presence of fully nonlinear relaxation terms.  

The outline of this paper is as follows. In the rest of this introductory
section, we review several model problems with linear and 
nonlinear relaxation, including Kawashima-LeFloch's model \cite{KL}. 
In Section~2, we  classify IMEX R--K schemes presented in literature
and indicate 
some limitations of these schemes when applied to hyperbolic
systems with diffusive relaxation. 
Section~3 is devoted to an introduction of partitioned and additive Runge-Kutta schemes for 
linear relaxation models \cite{BPR,BR}. In Sections~4 and 5, a particular type of IMEX R--K scheme, called
AGSA(3,4,2) (and introduced first in \cite{BR}), is analyzed and applied to the Euler equations with
linear friction and to a model coupling the Euler equations with a 
radiative transfert equation (investigated earlier in \cite{BLT} with a different numerical method). In Section~6, we describe and analyze the nonlinear relaxation model from \cite{KL}. 
Concluding remarks are finally presented in Section~7. 

\subsection{Models of interest}

Relaxation effects are an important feature of models arising, for instance, in compressible fluid dynamics and, specifically, we are interested here in nonlinear hyperbolic systems with fully nonlinear relaxation terms in the limit of, both, infinitely stiff relaxation and arbitrarily late time. In this limit, the dynamics of the flow is governed by effective systems of parabolic type, whose diffusion might also be degenerate and/or fully nonlinear. Solving such equations numerically is very challenging due to the stiffness of the problem in, both, the convection and the relaxation parts. This problem is currently under very active study and the main challenge is to design schemes that allow for realistic stability conditions on the time--step while ensuring robustness and high--order accuracy. For background as well as recent material, we refer to \cite{BLT,BR,JPT, JP, PLF, NP}. 

The simplest example of interest, which we propose to refer to as the {\bf linear relaxation model,}  reads 
\be
\label{eq:relax1}
\aligned
 u_t + v_x & = 0,
\\
\eps^2 \, v_t + b(u)_x & = -v + q(u),
\endaligned
\ee
in which $u=u(t,x)$ and $v=v(t,x)$ are the unknown functions, and $b=b(u)$ and $q=q(u)$ are prescribed with $b'(u)>0$. This is a nonlinear hyperbolic system with two (distinct and real) characteristic speeds, that is, $\pm \sqrt{b'(u)}/\eps$. In the stiff relaxation limit $\eps  \to 0$, the characteristic speeds approach infinity and the behavior of solutions to (\ref{eq:relax1}) is (formally, at least) governed by the effective system 
\be
\label{eq:redproblem}
\aligned
u_t  + q(u)_x &=  b(u)_{xx},
\\
v &= q(u) - b(u)_x,
\endaligned
\ee
which can be derived by a Chapman--Enskog expansion in the parameter $\eps$ and turns out to be a noninear parabolic system ---since the diffusion coefficient $b'(u)>0$ in $(b'u) u_x)_x$ is positive. Observe that the so-called sub--characteristic condition \cite{Whitham} for system (\ref{eq:redproblem}) reads $\eps \, |q'(u)| < (b'(u))^{1/2}$ and puts a restriction on the speed $q'(u)$ of the effective equation with respect to the speeds $\pm \sqrt{b'(u)}/\eps$ of the original system: clearly, for this model, the sub--characteristic condition is automatically satisfied (for $\eps$ sufficiently small).

Fully nonlinear relaxation terms also arise (for instance in presence of nonlinear friction) and, therefore, in this work  
we also numerically study the following {\bf nonlinear relaxation model}, first introduced by Kawashima and LeFloch \cite{KL}, 
\be
\label{eq:relax2}
\aligned
 u_t + v_x & = 0,
\\
\eps^2 \, v_t + b(u)_x & = - |v|^{m-1} \, v + q(u),
\endaligned
\ee
where $m>0$ is a real parameter. This, again, is a strictly hyperbolic system of balance laws, provided $b'(u)>0$. The effective system associated with this model is found to be (with $\alpha = -1+1/m$)
\be
\label{eq:redproblem2}
\aligned
u_t  &=  \big( |- q(u) + b(u)_x |^\alpha(- q(u) + b(u)_x ) \big)_x,
\\
 |v|^{m-1} \, v &= q(u) - b(u)_x,
\endaligned
\ee 
which is a {\sl fully nonlinear parabolic} equation in $u$. It is natural to distinguish between three rather distinct behaviors: 
\be
\label{values-m}
\aligned
\text{sub--linear : } &  0 < m < 1, 
\\ 
\text{linear : } &  m = 1, 
\\ 
\text{super--linear : } & m>1.  
\endaligned
\ee
The equation \eqref{eq:redproblem2}, although parabolic in nature, need not regularize initially discontinuous initial data, and we may expect jump singularities (in the first--order derivative $u_x$) to remain in the late--time limit. 

Further models arising in compressible fluid dynamics will also be described later, but \eqref{eq:relax1} and \eqref{eq:relax2} provide us two prototypes of particular importance in order to develop efficient numerical tools for relaxation problems. 
\section{Implicit-Explicit (IMEX) Runge-Kutta schemes}
This section is devoted to a brief introduction on IMEX-Runge-Kutta
schemes. More detailed description will be found in several papers on
the topic, such as \cite{ARS,CK,PaRu}.

Let us consider a Cauchy problem for a system of ordinary differential
equations of the form 
\begin{eqnarray}\label{odeIMEX}
   y' = f(y) + g(y), \quad y(0) = y_0, \quad t\in[0,T]
\end{eqnarray}
with $y(t)\in\mathbb{R}^N$, and $f$, $g$ two Lipschitz continuous functions
$\mathbb{R}^N\to \mathbb{R}^N$. 

Since Runge-Kutta schemes are one-step methods, it is enough to show
how to compute the solution after one time step, say from time $0$ to
time $\Delta t = h$.  An Implicit-Explicit (MEX) Runge-Kutta scheme has the
form 
\begin{eqnarray*}
 Y_i & = & y_0 + h \sum_{j=1}^{i-1} \tilde{a}_{ij} f(t_0+\tilde{c}_jh, Y_j) +
                 h \sum_{j=1}^i a_{ij} \frac{1}{\varepsilon} g(t_0+c_jh, Y_j),
                                                                 \label{eq:RKEI1}\\
 y_1 & = & y_0 + h \sum_{i=1}^{s} \tilde{b}_{i} f(t_0+\tilde{c}_ih, Y_i)  +
                 h \sum_{i=1}^s b_{i} \frac{1}{\varepsilon} g(t_0+c_ih, Y_i).
                                                                 \label{eq:RKEI2}
\end{eqnarray*}
where
$\tilde{A} = (\tilde{a}_{ij})$, $\tilde{a}_{ij}=0$, 
$j \geq i$ and $A = (a_{ij})$:
$s \times s$ are lower triangular matrices,
while $\tilde{c}, \tilde{b}, c, b \in \mathbb{R}^s$ are $s$-dimensional vectors. 
IMEX-RK schemes schemes are usually represented by a Double Butcher {\em tableau:}
\[
\begin{array}{c|c}
              \tilde{c} & \tilde{A} \\
              \hline \\
              & \tilde{b}^T
\end{array} \qquad
\begin{array}{c|c}
              c & A \\
              \hline \\
              & b^T
\end{array}
\]

The fact that the implicit scheme is diagonally implicit (DIRK) 
makes the implementation of IMEX simpler, and ensures $f$ is
effectively explicitly computed.
Order conditions for IMEX schemes can be derived by matching Taylor expansion of exact and numerical solution up to a
given order $p$, just as in the case of standard Runge-Kutta schemes.  

In addition to standard order conditions on the two schemes,
additional \emph{coupling conditions} may appear, whose number
increases very quickly with the order of the scheme. 
However, if $\tilde{c} = c$ and $\tilde{b} = b$, then there are no additional
conditions for IMEX-RK up to order $p=3$. For more details  on the
order conditions consult \cite{CK}. 

\subsection{Classification of IMEX R--K schemes}
IMEX R--K schemes present in the literature can be classified in three different types characterized by the structure of the matrix $A = (a_{ij})_{i,j=1}^s$ of the implicit scheme. Following \cite{Bosc-07}, we will rely on the following notions. 

\begin{definition}\label{Def:A}
An IMEX R--K method is said to be {\bf of type A} (cf.~\cite{PaRu}) if the matrix $A \in R^{s \times s}$ is invertible.
It is said to be {\bf of type CK} (cf.~\cite{CK}) if the matrix $A \in R^{s \times s}$ can be written in the form 
\begin{eqnarray*}
A = \left(\begin{array}{ll} 0 & 0\\
                        a  &  \widehat{A}\end{array}\right), 
               \end{eqnarray*}         
in which the matrix $\widehat{A} \in R^{(s-1) \times (s-1)}$ invertible.
Finally, it is said to be {\bf of type ARS} (cf.~\cite{ARS}) if it is a special case of the type CK with the vector $a = 0$.
\end{definition} 

The schemes CK are very attractive since they allow some simplifying assumptions, that make order conditions easier to treat, therefore permitting the construction of higher order IMEX R--K schemes. On the other hand, schemes of type A are more amenable to a theoretical analysis, since the matrix A of the implicit scheme is invertible. 

\subsection{A Simple Example of first IMEX R--K schemes}

The simplest IMEX-RK are obtained combining explicit and implicit 
Euler method. 
Two versions are possible applied to system (\ref{odeIMEX}), namely
\begin{eqnarray}\label{Afirst}
\mbox{SP(1,1,1), \cite{PaRu}:}
  \begin{array}{l}
    Y_1  =  y_0 + h \, g(Y_1)\\
    y_1  \> =  y_0 + h \, (f(y_0) + g(Y_1)),
  \end{array}
\end{eqnarray}
\begin{eqnarray}\label{ARSfirst}
\mbox{ARS(1,1,1), \cite{ARS}:}
  \begin{array}{rcl}
      y_1 & = & y_0 + h\, (f(y_0) + g(y_1)).
  \end{array}
  \end{eqnarray}
The corresponding Butcher tables are reported below: 
\[
\mbox{SP(1,1,1)}: \quad
     \begin{array}{c|c}
      0   & 0 \\
      \hline
          & 1
    \end{array}\quad
    \begin{array}{c|c}
      1 & 1        \\
      \hline
        & 1
    \end{array} \qquad
\mbox{ARS(1,1,1)}:\quad
     \begin{array}{c|cc}
      0   & 0 & 0 \\
      1   & 1 & 0 \\
      \hline
          & 1 & 0
    \end{array}\quad
    \begin{array}{c|cc}
      0   & 0 & 0  \\
      1   & 0 & 1 \\
      \hline
          & 0 & 1
    \end{array}
\]

\subsection{Limitations of IMEX Runge--Kutta schemes} 

In general, Implicit--Explicit (IMEX) Runge--Kutta (R--K) sche\-mes \cite{ARS, Bosc-07,Bosc-09, CK, PaRu} represent a powerful tool for the time
discretization of stiff systems. However, in the hyperbolic--to--parabolic limit under consideration here, the characteristic speeds of the hyperbolic part is of order $1/\eps$ and standard IMEX Runge--Kutta schemes developed for hyperbolic problems with stiff relaxation
\cite{PaRu, Bosc-10} are not very efficient, since the standard CFL condition requires the relation $\Delta t=\mathcal{O}(\eps \Delta x)$ between the time mesh--size and the space mesh--size. Clearly, in the diffusive regime $\eps << \Delta x$, and we are led to a too restrictive condition, since it is expected that, for an explicit method, a parabolic--type stability condition $\Delta t =\mathcal{O}(\Delta x^2)$ should suffice.

Most existing works \cite{JP,JPT,NP} on asymptotic preserving schemes for nonlinear hyperbolic system with diffusive relaxation are based on a (micro--macro) {\sl decomposition of the hyperbolic flux into stiff and non--stiff components}: the stiff hyperbolic component is added to the relaxation term, which is treated {\sl implicitly;} in the limit of infinite stiffness, such schemes become consistent and {\sl explicit} schemes for the diffusive limit equation \cite{JPT, Klar, LS, LM, NP} and, therefore, suffer from the usual stability restriction $\Delta t =\mathcal{O}(\Delta x^2)$. 

Schemes that do avoid the parabolic time--step restriction and provide {\sl fully implicit} solvers (in the case of transport equations) have been analyzed in \cite{BPR}, where a new formulation of the problem (\ref{eq:relax1}) in the case $b(u) = u$ and $q(u) = 0$, 
was introduced, which was based on the addition of {\sl two opposite diffusive terms} in the stiff limit. One term is added to the hyperbolic component and makes it non--stiff (therefore allowing for an explicit treatment), while the other term must be treated implicitly. The remaining component of the hyperbolic term is formally treated implicitly. The resulting scheme is consistent and relaxes to an implicit scheme for underlying diffusion limit, therefore avoiding the typical parabolic restriction of previous methods.

The drawback of the above strategy is that its requires an {\sl implicit} treatment of some components of the hyperbolic flux. Even if such implicit term can be computed by closed formula in most cases, this procedure requires a reformulation of the discretization of the whole system. It is thus desirable to construct IMEX R--K schemes in which the whole hyperbolic part is treated {\sl explicitly.} At a first sight this seems an impossible task since the characteristic speeds are unbounded. Yet, in \cite{BR}, the authors showed that, provided suitable conditions are imposed on the IMEX R--K coefficients, then one can indeed construct schemes which are
\begin{itemize}

\item  fully explicit in the hyperbolic part, and 

\item converge to an implicit high--order scheme in the diffusion limit. 

\end{itemize}
This is this direction that we follow in the present work and further develop in the context of, both, the linear relaxation model \eqref{eq:relax1}, the fully nonlinear model \eqref{eq:relax2}, and several models arising in fluid dynamics. We exhibit new decompositions that lead us to high--order asymptotic preserving methods covering the hyperbolic--to--parabolic regime under consideration. A key ingredient is to keep the implicit feature of the scheme with respect to certain components that ``degenerate''; for instance, \eqref{eq:relax2} when $m>1$ has a degenerate diffusion when $- q(u) + b(u)_x$ is small or even vanishes.


\section{Partitioned and Additive Runge--Kutta schemes for the linear relaxation model}

In developing new IMEX R--K schemes in \cite{BR,BPR}, the system ({\ref{eq:relax1}}) was considered and 
written in the form
\begin{eqnarray}\label{refsystem}
  u_t &=& -v_x ,                  \nonumber \acapo
              &  &                    \acapo
 \varepsilon ^2 v_t &=&  -( b(u)_x + v-q(u )).      \nonumber
\end{eqnarray}
A very simple scheme for the numerical treatment of such a problem is
based on the observation that the first equation does not contain any
stiffness, so the first variable can be updated in time, and the
new value can be replaced in the second equation. 

By using a method of lines approach (MOL), we discretize system (\ref{refsystem}) in space by a uniform mesh $\left\{x_i\right\}_{i  =1}^N$  and  $U_i(t) \approx u(x_i, t)$, $V_i(t) \approx v(x_i, t)$. We obtain a large system of ODE's 
\[
   U_t = F(U,V), \quad \varepsilon^2 V_t = G(U) - V,
\]
where, in this case, $F(U,V) = -\mathcal{D}_h V$, $G(U) = q(U) - \mathcal{D}_hb(U)$, $\mathcal{D}_h$
denoting a discrete space derivative operator. 
In the limit $\varepsilon\to 0$ the system becomes 
\[
   U_t = \hat{F}(U), \quad \mbox{where} \quad \hat{F}(U) = F(U,G(U)).
\]
The IMEX Euler scheme (/ref{Afirst}) applied to the system
becomes
\begin{eqnarray}\label{systemFirst}
\begin{array}{l}
   U^{n+1} =  U^n + \Delta t F(U^n,V^n), \\
   V^{n+1}  =   V^n + \Delta t G(U^{n+1}) - \Delta t V^{n+1}, 
   \end{array}
\end{eqnarray}
and the second equation can be explicitly solved for $V^{n+1}$
\[
   V^{n+1} = \frac{\varepsilon^2V^n + \Delta t G(U^{n+1})}{\varepsilon^2+\Delta t}.
\]
where we have discretized the interval of integration by a time mesh $\left\{t_n\right\}_{n  =1}^\mathcal{N}$ and $U^n \approx  U(t^n)$.
As $\varepsilon \rightarrow 0$, $V^{n+1} = G(U^{n+1})$ and therefore $U^{n+1} = U^{n}+\Delta t \hat{\hat{F}}(U^{n})$, i.e. the scheme automatically become an explicit Euler
method 
applied to the space discretized limit equation. 
Note that in system (\ref{systemFirst}), since the variable $U$ is updated in the first equation,
it is explicitly computed, although it is formally implicit. 

In the case $b(u) = u$ and $q(u)=0$, the method would relax to explicit Euler scheme applied to the diffusion equation, thus suffering the usual parabolic 
CFL stability condition $\Delta t \leqslant \Delta x^2/2$.
This approach will be denoted as \emph{partitioned approach} and we will denote the corresponding IMEX R--K schemes, {\bf IMEX-I} R--K schemes .  

The simplest scheme outlined above is only first order accurate in
time. However, higher order extensions are possible, by adopting IMEX
type time discretization.  
In the previous approach there may be a subtle difficulty when it
comes to applications, namely it is not clear how to identify the
hyperbolic part of the system, i.e.\ what is the term that should be
included in the numerical flux. 

For practical applications, it would be very nice to treat the whole
term containing the flux explicitly, while reserving the implicit
treatment only to the source, according the following scheme:
\[
\begin{array}{r}
u_t  \\
v_t  \\
{~}
\end{array} 
\begin{array}{c}
 = \\
 = \\
{~}
\end{array}
\begin{array}{c}
  \framebox{$
    \begin{array}{c} - v_x\\ -u_x/\varepsilon^2 
    \end{array}$} \\
  {\rm [Explicit]}
\end{array}
\hspace{2mm}
\begin{array}{c}
  \framebox{$
    \begin{array}{cl}
      &   \\
      - & v/\varepsilon^2 
    \end{array}
    $} \\
  {\rm [Implicit]}
\end{array}
\begin{array}{c}
\\[-2mm]
 \quad {\rm (Additive)}\\[2mm]
\\
\end{array}
\]
We call such an approach {\em additive} and the corresponding schemes
are denoted {\bf IMEX-E}.  Such schemes should be easier to apply, since the fluxes retain
their original interpretation. However, this approach seems to be hopeless, due to diverging speeds.

Let us consider the simple relaxation system with $b(u) = u$,
and $q=0$. By applying the Implicit-Explicit Euler scheme ARS(1,1,1), \rf{ARSfirst} 
we obtain 
\begin{eqnarray*}
   u^{n+1} & = & u^n + \Delta t v_x^n,
 \\
   v^{n+1} & = & v^n - \Delta t \, u_x^n/\varepsilon^2 - \Delta t \, v^{n+1}/\varepsilon^2.
\end{eqnarray*}
Now, as $\varepsilon \to 0$, one  has $v^{n+1}\to -u_x^n$. 
Replacing this expression in the first equation gives
\[
   u^{n+1} = u^n + \Delta t u_{xx}^{n-1}, 
\]
which provides a {\em consistent discretization} of the limit
equation.

Performing Fourier stability analysis on this scheme, one obtains a
stability condition that depends on $\varepsilon$.
As $\varepsilon \to 0$ one obtains
\[
   \begin{array}{c}{\Delta t \xi^2\leq 1} \\[0.3cm]
       \mbox{continuous in space}
    \end{array}
    \quad 
      \mbox{or} \quad
   \begin{array}{c} \displaystyle
   {{\Delta t}\leq {\Delta x^2}} \\[0.3cm]
    \mbox{central difference}.
    \end{array}
\]
Unfortunately, when trying this approach with other IMEX schemes
available in the literature, even with the first order SP(1,1,1), one
observes a lack of consistency (for details, see\cite{BR}). 

In this paper, the authors concentrated on developing IMEX R--K schemes of type A, since 
 they are easier to analyze with respect to the other types. They started the analysis by introducing a property which is important in order 
to guarantee the asymptotic preserving property.  IMEX R--K schemes that satisfy this property are called  {\em globally stiffly
  accurate}, that is the last row of matrices $\tilde{A}$ and $A$ are,
respectively, equal to $\tilde{b}^T$ and $b^T$. 

\begin{definition} \label{1def}
An IMEX R--K scheme 
is said to be  {\bf globally stiffly accurate} if 
$$
\bt = e_s^T A, 
\qquad 
\tbt = e_s^T\tA
\qquad \text{ with }
e_s = (0,...,0,1)^T, \quad
c_s = \widetilde{c}_s = 1, 
$$
i.e.~the numerical solution coincides with the solution obtained in the last internal step of the Runge--Kutta scheme. 
\end{definition}

Note that this condition is satisfied by ARS(1,1,1), (\ref{ARSfirst}), but
not by SP(1,1,1), (\ref{Afirst}).
An IMEX-E R--K scheme applied to this approach   
will relaxes to an explicit scheme for the limit parabolic equation, 
thus suffering from the classical parabolic CFL stability restriction
on the time step.


\subsection{Removing parabolic stiffness}

The schemes constructed with the approachs outlined above will converge
to an explicit scheme for the limit diffusion equation, and therefore
they are subject to  the classical parabolic CFL restriction 
$\Delta t \leq C \Delta x^2$. In order to overcome such a restriction we adopt a
penalization technique, based on adding two opposite terms to the
first equation, and treating one explicitly and one implicitly.

Let us consider system (\ref{eq:relax1}) and adding and subtracting the same term $\mu(\varepsilon) b(u)_{xx}$ on the right hand side we 
obtain:
\be
\label{startPoint}
\aligned
   u_t   &=&  - (v +  \mu(\varepsilon) b(u)_{x})_x + \mu(\varepsilon) b(u)_{xx}, 
   \\
   v_t   &=&  -  \frac{1}{\varepsilon^2} (b(u)_x + v - q(u)).
\endaligned
\ee
Here $\mu = \mu(\varepsilon) \in [0, 1]$ is a free parameter such that $\mu(0) = 1$. The idea is that, since the quantity $v + b(u)_x$ is close to $q(u)$ as $\varepsilon \to 0$, the first term on the right hand side can be treated explicitly in the first equation, while the term $p(u)_xx$ will be treated implicitly. This can be done naturally by using an IMEX R--K approach.

However, if we want the scheme to be accurate also in the
cases in which $\varepsilon$ is not too small, then we need to add two main
ingredients:
\begin{itemize}
\item[1.] no term should be added when not needed, i.e.\ for large enough
  values of $\varepsilon$, since in such cases the additional terms
  degrade the accuracy; this
  can be achieved by letting $\mu(\varepsilon)$ decrease as $\varepsilon$
  increases. A possible choice, for example, is 
\[
  \mu(\varepsilon) = \exp(-\varepsilon^2/\Delta x).
\]
\item[2.] when the stabilizing effect of the dissipation vanishes, i.e. as
  $\mu \to 0$, then central differencing is no longer suitable, and one
  should adopt some upwinding; this can be obtained for example by
  blending central differencing and upwind differencing as 
\[
   D_x = (1-\mu)D_x^{\rm upw} + \mu D_x^{\rm cen}.
\]
\end{itemize}

Now from a numerical point of view, we may treat system (\ref{startPoint}) in two different approaches according to whether the hyperbolic term $b(u)_x$ in the second equation is treated implicitly or explicitly.
The first one, called the \emph{BPR approach} (For more details as well as some rigorous analysis, see \cite{BPR}.) is 
\begin{eqnarray}\label{systempq2} 
  u_t          & = &\boxed{- (v+\mu(\eps) b(u)_x)_x}_{\text{ Explicit}} + \boxed{\mu(\eps) b(u)_{xx}}_{\text{ Implicit}}   \nonumber        \acapo
               &   &                        \acapo
 \eps^2 v_t & = & \boxed{-(b(u)_x +v - q(u))}_{\text{ Implicit}}        \nonumber
\end{eqnarray}
and the corresponding IMEX R--K schemes are {\bf IMEX-I} R--K ones in order to remind that the term containing $b(u)_x$ in the
second equation is implicit, in the sense that it appears at the new
time level. 
Observe that the term $b(u)_x + v - q(u)$ appearing in the second equation 
is formally treated implicitly, but in
practice the term $b(u)_x$ is computed at the new time from the first equation, so
it can in fact explicitly computed. 

The second approach (cf.~\cite{BR} for further details) is given by 
\begin{eqnarray}\label{systempq} 
  u_t          & = &\boxed{- (v+\mu(\eps) b(u)_x)_x}_{\text{ Explicit}} + \boxed{\mu(\eps) b(u)_{xx}}_{\text{ Implicit}}   \nonumber        \acapo
               &   &                        \acapo
 \eps^2 v_t & = & \boxed{-b(u)_x}_{\text{ Explicit}} -\boxed{(v - q(u))}_{\text{ Implicit}}, 
     \nonumber
\end{eqnarray}
which is refered to as the  \emph{BR} approach; he corresponding schemes are {\bf IMEX-E}.

Observe that  when $q(u) \neq 0$, the system relaxes to a convection--diffusion equation. The terms on the left hand side are treated explicitly, while the terms on the right--hand side are treated implicitly.  In fact as $\eps \to 0$, the scheme becomes an IMEX R--K scheme for the limit convection-diffusion equation, in which the convection term is treated explicitly and the diffusion one is treated implicitly \cite{BPR,BR}. 


\subsection{Additional order conditions}

In \cite{BPR,BR} additional order conditions, called {\em algebraic conditions} are derived by studying the behavior of IMEX R--K schemes applied to the  relaxation system (\ref{startPoint}) when $\eps \to 0$, in particular when $\eps = 0$. These algebraic conditions guarantee the correct behavior of the numerical solution in the limit $ \varepsilon \to 0$ and maintain the accuracy in time of the scheme. We obtained such algebraic order conditions using the same technique adopted in \cite{Bosc-09}, i.e.~ by comparing the Taylor series of the numerical solution with that of the exact solution. 

In \cite{BPR} these order conditions are derived for {\bf IMEX-I} R--K schemes and in \cite{BR} for {\bf IMEX-E} R--K ones.
Due to the number of these order conditions, it is necessary to use
several internal levels for IMEX R--K schemes in order to obtain a given order. Several results and a rigorous analysis about that can be found in \cite{BPR, BR}. 
Most numerical tests are reported in \cite{BPR} for the BPR approach and in \cite{BR} for the BR approach and 
the results are compared with those obtained by other methods
available in the literature. 

Here we report, as an example, the application of the second order {\bf IMEX-E} R--K scheme developed in \cite{BR} and called AGSA(3,4,2). We note that this scheme satisfy all the algebraic order conditions derived in \cite{BR}.  Now we perform several numerical tests, considered earlier in \cite{BLT} with a different 
numerical strategy.

\section{Euler equation with linear friction} 

The Euler model with friction reads
\begin{eqnarray}\label{Euler}
\begin{array}{ll}
\displaystyle \eps \rho_t + (\rho w)_x = 0,\\
\displaystyle \eps (\rho w)_t + (\rho w^2 + p(\rho))_x = -\frac{1}{\eps}\rho w,
\end{array}
\end{eqnarray}
where $\rho>0$ denotes the density and $w \in \mathbb{R}$ the velocity of a compressible fluid. The pressure function $p:\mathbb{R}^+ \to \mathbb{R}^+$ is assumed sufficiently regular and satisfies $p'(\rho)>0$ so that the first-order homogeneus system associated with (\ref{Euler}) is strictly hyperbolic.

Now using the new variable $v = w/\eps$ the system becomes
\begin{eqnarray}\label{Euler2}
\begin{array}{ll}
\displaystyle \rho_t + (\rho v)_x = 0,\\
\displaystyle (\rho v)_t + (\rho v^2 + \frac{p(\rho)}{\eps^2})_x = -\frac{1}{\eps^2}\rho v.
\end{array}
\end{eqnarray}
In this case the diffusive regime associated to the Euler model with friction is described when $\eps \to 0$ by the equation 
\begin{eqnarray}\label{limitequation}
\rho_t = (p(\rho))_{xx}.
\end{eqnarray}

Now, from a numerical point of view in order to test our examples we adopt a technique similar to the one illustrated
in \cite{BPR,BR}. Indeed, the more classical techniques for instance proposed in \cite{PN, JPT, BLT} provide in the asympotic diffusive limit a scheme that converges to an explicit scheme for the parabolic equation. But such schemes suffer from the
standard CFL restriction $\Delta t = \mathcal{O}(\Delta x^2)$ when $\eps \to 0$. Then, in order to remove such restriction as suggested  in \cite{BPR,BR}, this technique consists in adding and subtracting the same term to the first equation
of system (\ref{Euler2}):
\begin{eqnarray}\label{Euler2bis}
\begin{array}{ll}
\displaystyle \rho_t =  - \big( \rho v +\mu(\eps) (p(\rho)_x \big)_{x} -  \mu(\eps) \, p(\rho)_{xx},\\
\displaystyle (\rho v)_t  =- \Big( \rho v^2 + \frac{p(\rho)}{\eps^2} \Big)_x -\frac{1}{\eps^2}\rho v.
\end{array}
\end{eqnarray}
Here $\mu$ is such that $\mu(\eps) : \mathbb{R}^+ \to [0, 1]$ and $\mu(0) = 1$. When $\eps$ is not small there is no
reason to add and subtract the term $(p(\rho))_{xx}$, therefore $\mu(\eps)$ will be small in such a
regime, i.e. $\mu(\eps) = 0$. A possible choice for $\mu(\eps)$ can be
\begin{eqnarray}
\mu(\eps) =\begin{cases} 
1, \ & \eps <\Delta x,
\\
0, \quad  & \eps \geq \Delta x.
\end{cases}
\end{eqnarray}
With this approach, the idea is that as $\eps \to 0$, the quantity $(\rho v) +\mu(\eps) p(\rho)_x \to 0$. Therefore
such a term can be treated explicitly, while the other term $p(\rho)_{xx}$ will be treated
implicitly in the first equation, i.e. we will treat the hyperbolic part of the system
$(((\rho v) +\mu(\eps) p(\rho)_x)_{x}, (\rho v^2 + \frac{p(\rho)}{\eps^2})_x)^T$ explicitly, while the relaxation one $((p(\rho))_{xx},-\frac{1}{\eps^2}\rho v)^T$ implicitly, respectively.

We performed our scheme AGSA(3,4,2) introduced in the previous section for the numerical experiment considering a parabolic $\Delta t$, i.e. $\Delta t \propto \Delta x^2$ by solving system (\ref{Euler2}) and hyperbolic $\Delta t$, i.e. $\Delta t \propto \text{CFL} \, \Delta t$ by solving (\ref{Euler2bis}). We show only the last case since we obtained similar results with the first one. 

\begin{figure}\label{Fig:1}
\centering 
\includegraphics[width=0.7\textwidth]{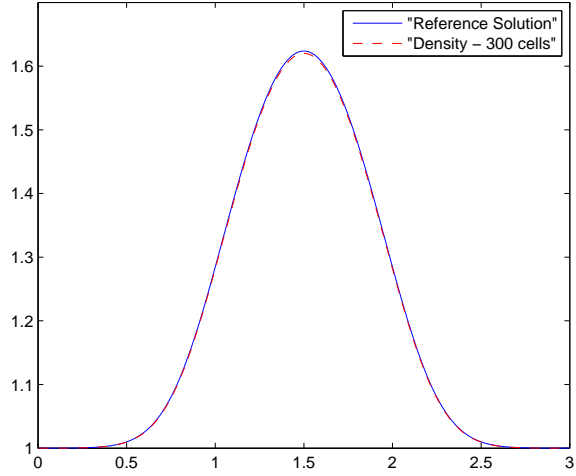}
\caption{Comparison between {\bf reference solution} (full line) obtained by the numerical solution or (\ref{limitequation}) of (4.3) in \cite{BLT} computed with 600 cells $\Delta x = 10^{-2}$and $\Delta t \propto \Delta x^2 \cong 5.10^{-5}$ versus {\bf proposed scheme solution} (dashed line) ASA(3,4,2) computed with 300 cells and $\Delta t = 0.1\Delta x \cong 1.10^{-3}$. The solutions are compared at time $T_f = 2.10^4\eps$ with $\eps = 10^{-3}$.}
\label{fig:1}
\end{figure}

We now solve the (\ref{Euler2}) with the inital data
\begin{eqnarray}\label{InDat}
(\rho,\rho v)^{T} =\begin{cases} 
(2, 0)^{T}, & x\in [1.2, 1.8],
\\
(1, 0)^{T}, & \textrm{otherwise}.
\end{cases}
\end{eqnarray}
Furthermore we choose the simple pressure law $p(\rho) = \rho^2$ and $\eps = 10^{-3}$. In Figure \ref{Fig:1}, the numerical solution on a $300$ cells grid with $\Delta x = 10^{-2}$ obtained by the proposed scheme, is compared at time $t_{final} = 2.10^4\eps$, with a numerical approximation of (\ref{limitequation}), (or see (4.3) in \cite{BLT}) computed with $600$ cells. 


\section{Coupled Euler--radiative transfert model} 

\begin{figure}\label{Fig:2}
\centering 
\includegraphics[width=0.7\textwidth]{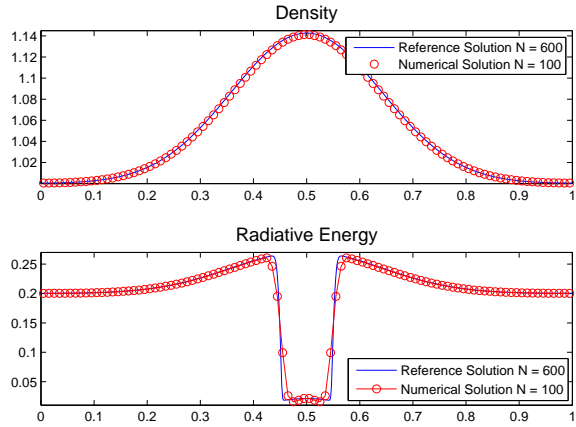}
\caption{Comparison between {\bf reference solution} (full line) obtained by (\ref{limitequation2}) or the numerical solution of (4.6) in \cite{BLT} computed with 600 cells, $\Delta x = 10^{-2}$ and $\Delta t \propto \Delta x^2 \cong 5.10^{-5}$ versus {\bf proposed scheme solution} ASA(3,4,2) computed with 100 cells and $\Delta t = 0.1\Delta x \cong 1.10^{-3}$. The solutions are compared at time $T_f = 0.029$ with $\eps = 10^{-3}$.}
\label{fig:2}
\end{figure}

The second example proposed here is a system that degenerates into a system of diffusion equation of dimension $n>1$ and in order to do this we couple the Euler model (\ref{Euler}) with the $M1$ model proposed in \cite{BLT} as follows
\begin{eqnarray}\label{EulerM1}
\begin{array}{ll}
\displaystyle \eps \rho_t + (\rho w)_x = 0,\\
\displaystyle \eps (\rho w)_t + (\rho w^2 + p(\rho))_x = -\frac{\kappa}{\eps}\rho w + \frac{\sigma}{\eps} \bar{f},\\
\displaystyle \eps e_t + (\rho \bar{f})_x = 0,\\
\displaystyle \eps \bar{f}_t + (\chi\left(\frac{\bar{f}}{e}\right)e)_x = -\frac{\sigma}{\eps}\bar{f}.
\end{array}
\end{eqnarray}
 Now introducing the new variables $v = w/\eps$ and $f = \bar{f}/\eps$ we get the system
\begin{eqnarray}\label{EulerM12}
\begin{array}{ll}
\displaystyle \rho_t + (\rho v)_x = 0,\\
\displaystyle (\rho v)_t + (\rho v^2 + \frac{p(\rho)}{\eps^2})_x = -\frac{\kappa}{\eps^2}\rho v + \frac{\sigma}{\eps^2} f,\\
\displaystyle  e_t + (\rho f)_x = 0,\\
\displaystyle f_t + \frac{1}{\eps^2}(\chi\left(\frac{\eps f}{e}\right)e)_x = -\frac{\sigma}{\eps^2}f.
\end{array}
\end{eqnarray}
Here $\kappa$ and $\sigma$ denote positive constants. Furthermore the following pressure law: $p(\rho) = C_p \rho^{\eta}$ with $C_p\ll 1$ and $\eta > 1$.
The asymptotic diffusive regime (i.e. $\eps \to 0$) of the system (\ref{EulerM12}) is given by
\be
\label{limitequation2}
\aligned
 \rho_t &= \frac{p(\rho)_{xx}}{\kappa} + \frac{e_{xx}}{3\kappa},
 \\
e_t &= \frac{e_{xx}}{3 \sigma}.
\endaligned
\ee
Now similarly as made for the previous example here in order to remove time step parabolic restriction we adding and subtract the same term to the first and third equation
\begin{eqnarray}\label{EulerM12bis}
\begin{array}{ll}
\displaystyle \rho_t = - \left((\rho v) + \frac{1}{\kappa}\left( p(\rho)_x + \frac{1}{3} e_x\right)\right)_x + \frac{1}{\kappa}\left( p(\rho)_{xx} + \frac{1}{3} e_{xx}\right),
\\
\displaystyle (\rho v)_t + (\rho v^2 + \frac{p(\rho)}{\eps^2})_x = -\frac{\kappa}{\eps^2}\rho v + \frac{\sigma}{\eps^2} f,\\
\displaystyle  e_t = - ((\rho f) + \frac{e_{x}}{3 \sigma})_x + \frac{e_{xx}}{3 \sigma},\\
\displaystyle f_t - \frac{1}{\eps^2}(\chi\left(\frac{\eps f}{e}\right)e)_x = -\frac{\sigma}{\eps^2}f,\\
\end{array}
\end{eqnarray}
Similarly for this system (\ref{EulerM12bis}) we perform our scheme considering the hyperbolic $\Delta t$, i.e. $\Delta t \propto \text{CFL} \, \Delta t$. 
We choose initial data given by
\begin{eqnarray}\label{InDat2}
\rho = 0.2, \ \ v = 0, \ \ f=0, \ \ e =\begin{cases} 
1.5, & x\in [0.45, 0.55],
\\
1, & \textrm{otherwise}.
\end{cases}
\end{eqnarray}
The parameters of the model are $\kappa = 2$, $\sigma=1$, $\eps = 10^{-3}$, $T_{final} = 0.029$, $C_p = 10^{-3}$ and $\eta = 2$. The numerical solution is computed with $100$ cells ($ \Delta x = 10^{-2}$) and compared to Figure \ref{Fig:2} with a reference solution obtained solving (\ref{limitequation2})  

In Figures~5.1 and 5.2, we observe that the second--order scheme AGSA(3,4,2) developed in \cite{BR} captures well the correct behavior of the solutions  in the diffusive regime where the numerical solutions match the reference solution very well.


\section{Kawashima-LeFloch's nonlinear relaxation model}

\subsection{Objective}

We are now in a position to tackle the nonlinear relaxation model (\ref{eq:relax2}) 
which is 
of central interest in the present paper. Without genuine loss of generality in the design of the numerical method, we assume that $b(u) = u$ and $q(u)= 0$ and, therefore, we study 
\be
\label{eq:relax2bis}
\aligned
 u_t + v_x & = 0,
\\
\eps^2 \, v_t + u_x & = - |v|^{m-1} \, v, 
\endaligned
\ee
together with the associated effective equation ($\varepsilon \to 0$)
\be
\label{eq:redproblem2bis}
u_t  =  \big( |u_x |^{\alpha} u_x  \big)_x, 
\qquad
\quad
|v|^{m-1} \, v = - u_x, 
\ee  
with  $\alpha = -1+1/m$. 

The order of convergence of the IMEX scheme applied to the semilinear system (\ref{eq:relax2bis}) for small $\varepsilon$ will be determined by comparing the numerical solutions to a (theoretical) {\em limiting solution} obtained by solving the effective equations (\ref{eq:redproblem2bis}) on a very fine mesh.

We are going to treat the nonlinear relaxation term {\sl implicitly,} which requires to solve a nonlinear algebraic equation. An efficient approach is obtained by solving it iteratively with the Newton method. Specifically, the nonlinear equation to be solved at each $s$--stage has the form 
\begin{eqnarray*}
F(V) = \varepsilon^2 \, V + C_1 \, V^m - \varepsilon^2 \, C_2  = 0,
\end{eqnarray*}
where $C_1$ and $C_2$ constants. As $\varepsilon \to 0$, we can neglect the term $V \varepsilon^2$ and, therefore, we propose the initial approximation 
\begin{eqnarray*}
V_0 = \Big( \frac{\varepsilon^2 C_2}{C_1} \Big)^{1/m},
\end{eqnarray*}
as a good choice in order to reduce the number of iterations. This is a relevant choice in both cases $m<1$ and $m \geq 1$.


\subsection{Computation of the limit solution}
\label{sec:reflim}

Equation (\ref{eq:redproblem2bis}) is a nonlinear diffusion equation. Let us introduce a very efficient method for the numerical solution of such an equation. The scheme we propose is stable, linearly implicit, and can be designed up to {\sl any order of accuracy.}  We propose to use a technique similar to the additive semi-implicit Runge-Kutta methods of Zhong \cite{Zhong}. The idea is to write equation as a system
\begin{eqnarray}
\aligned
\label{part}
u' = F(u,u), 
\endaligned
\end{eqnarray}
where the first argument is treated explicitly, while the second one is treated implicitly. Then the semi-implicit R--K method is implemented as follows. First compute the stage fluxes for $i = 1, \ldots, s$
\begin{eqnarray}
U^{\ast}_i = u^{n} + \Delta t \sum_{j =1}^{j-1} \tilde{a}_{ij} K_j, 
\qquad \quad
\bar{U}_i = u^{n} + \Delta t \sum_{j =1}^{j-1} a_{ij} K_j,
\end{eqnarray}
solve the equation
\begin{eqnarray*}
K_i = F(U^{\ast}_i, \bar{U}_i + \Delta t a_{ii}K_i), 
\end{eqnarray*}
and, finally update the numerical solution
\begin{eqnarray}\label{refsol2}
u^{n+1} = u^{n} + \Delta t \sum_{i =1}^{s} b_i K_i.
\end{eqnarray}

In our case $F(u^{\ast}, u)$ is given by  
$$
F(u^{\ast}, u) = (|u^{\ast}_x|^{\alpha}u_x)_x.
$$ 
It can be shown that such semi-implicit R--K scheme is indeed a particular case of IMEX R--K scheme \cite{PaRu2}. 

Observe that, for the spatial discretization of the limiting equation (\ref{eq:redproblem2bis}), we can use a centered scheme, such as 
\begin{eqnarray*}
\big( |u_x |^{\alpha} u_x ) \big)_x \approx 
 \Delta x \, 
\Bigg( |\delta U|^{\alpha}_{j+\frac{1}{2}}\frac{U_{j+1}(t)-U_j(t)}{\Delta x} - |\delta U|^{\alpha}_{j-\frac{1}{2}}\frac{U_{j}(t)-U_{j-1}(t)}{\Delta x}\Bigg),
\end{eqnarray*}
where $ |\delta U|^{\alpha}_{j+\frac{1}{2}} = \Big| \frac{U_{j+1}(t)-U_j(t)}{\Delta x} \Big|^{\alpha}$.
Observe that for the final numerical solution, we require, in particular, that $u^{\ast^{ n+1}} = u^{n+1}$, which is guaranteed by imposing the condition $b_i = \widetilde{b}_i$ for all $i$.

Here,  as an example, we propose to choose a classical second--order IMEX R--K scheme that satisfies the condition $b_i = \widetilde{b}_i$ for all $i$, e.g. the {\sl implicit--explicit midpoint method of type ARS} introduced in \cite{ARS} (to which we refer for further details). Applying such a method, we thus obtain the following (rather simple) scheme which requires only two function evaluations
\be
\label{stages_part}
\aligned
U^{\ast}_2 &= u_n + \frac{\Delta t}{2} F(u_n, u_n), 
\\
U_2 &= u_n + \frac{\Delta t}{2}F(U^{\ast}_2,U_2), 
\\
u_{n+1} &= 2 U_2 - u_n.
\endaligned
\ee 
Observe that $\mathcal{U}^{\ast}_1 = \mathcal{U}_1 = u_n$. 

Similar examples is possible to propose for other types of IMEX R--K schemes presented in the literature. Then, we consider the numerical solution (\ref{stages_part}) as a reference solution for the limiting equation, and we compare it with the numerical one computed from the system (\ref{eq:relax2bis}). 


\subsection{Comparison with the limit solution}
Here we apply first order scheme to system (\ref{eq:relax2bis}) with simple second order central differencing for the discretization of the space derivative. We choose the first order implicit-explicit Euler scheme of type ARS (\ref{ARSfirst}).  
For our numerical test we integrate the system up to the final time $T = 1$ and use periodic boundary conditions in $x \in [-\pi, \pi]$ with initial conditions $u(x,0) = \cos x$, $v(x,0) = \sin x$ and $\varepsilon = 10^{-4}$.
 
Numerical convergence rate is calculated by the formula
\begin{eqnarray*}
p = \log_2(E_{\Delta t_1}/E_{\Delta t_1}),
\end{eqnarray*} 
where $E_{\Delta t_1}$ and $E_{\Delta t_2}$ are the global errors computed with step $\Delta t_1 = \mathcal{O}(\Delta x_1^2)$ where $\Delta x_1 = 2 \pi / N$ and  $\Delta t_2 = \mathcal{O}(\Delta x_2^2)$ with $\Delta x_2 =  \pi /  N$.

\begin{figure}\label{Fig:3}
\centering 
\includegraphics[width=0.4\textwidth]{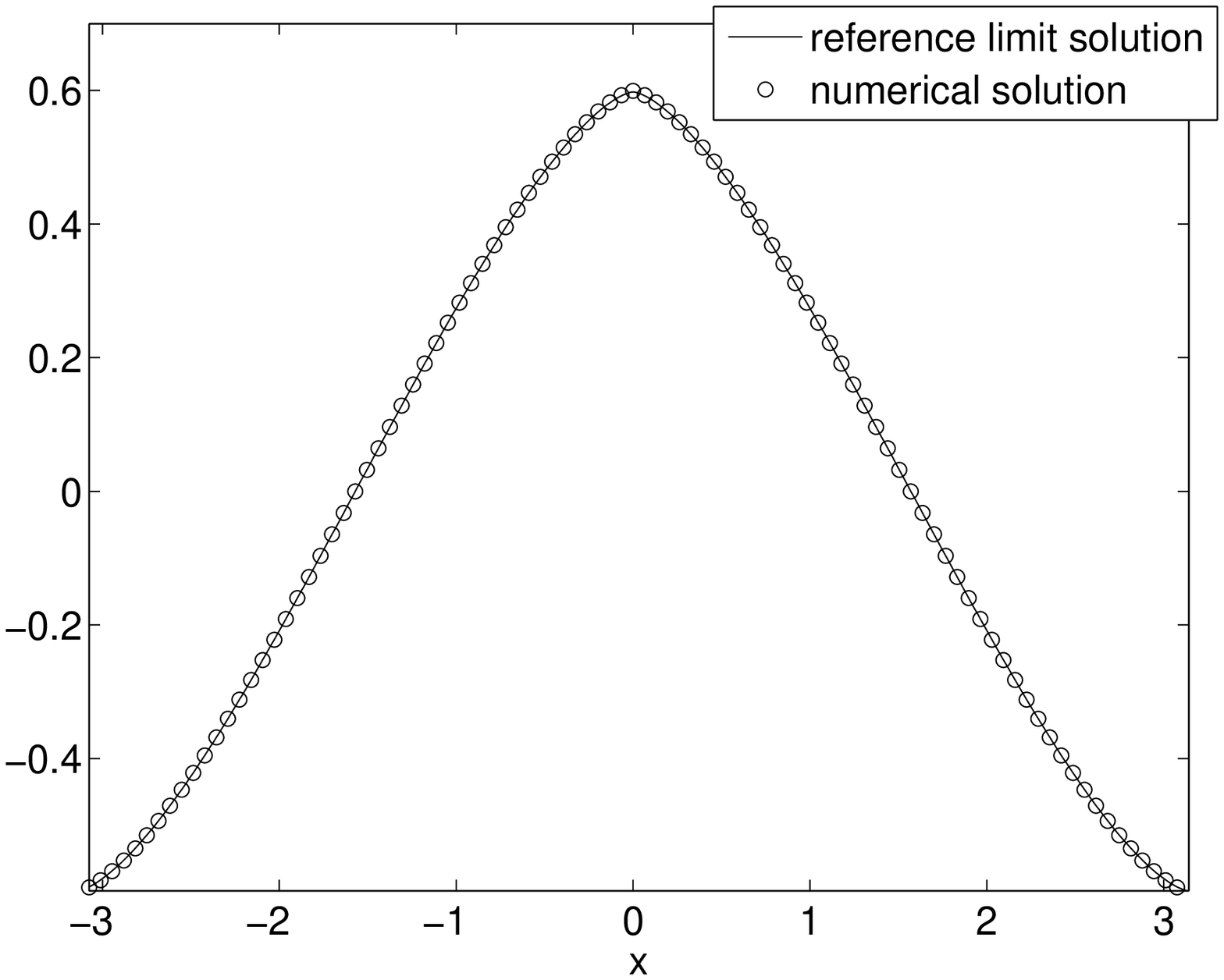}
\includegraphics[width=0.48\textwidth]{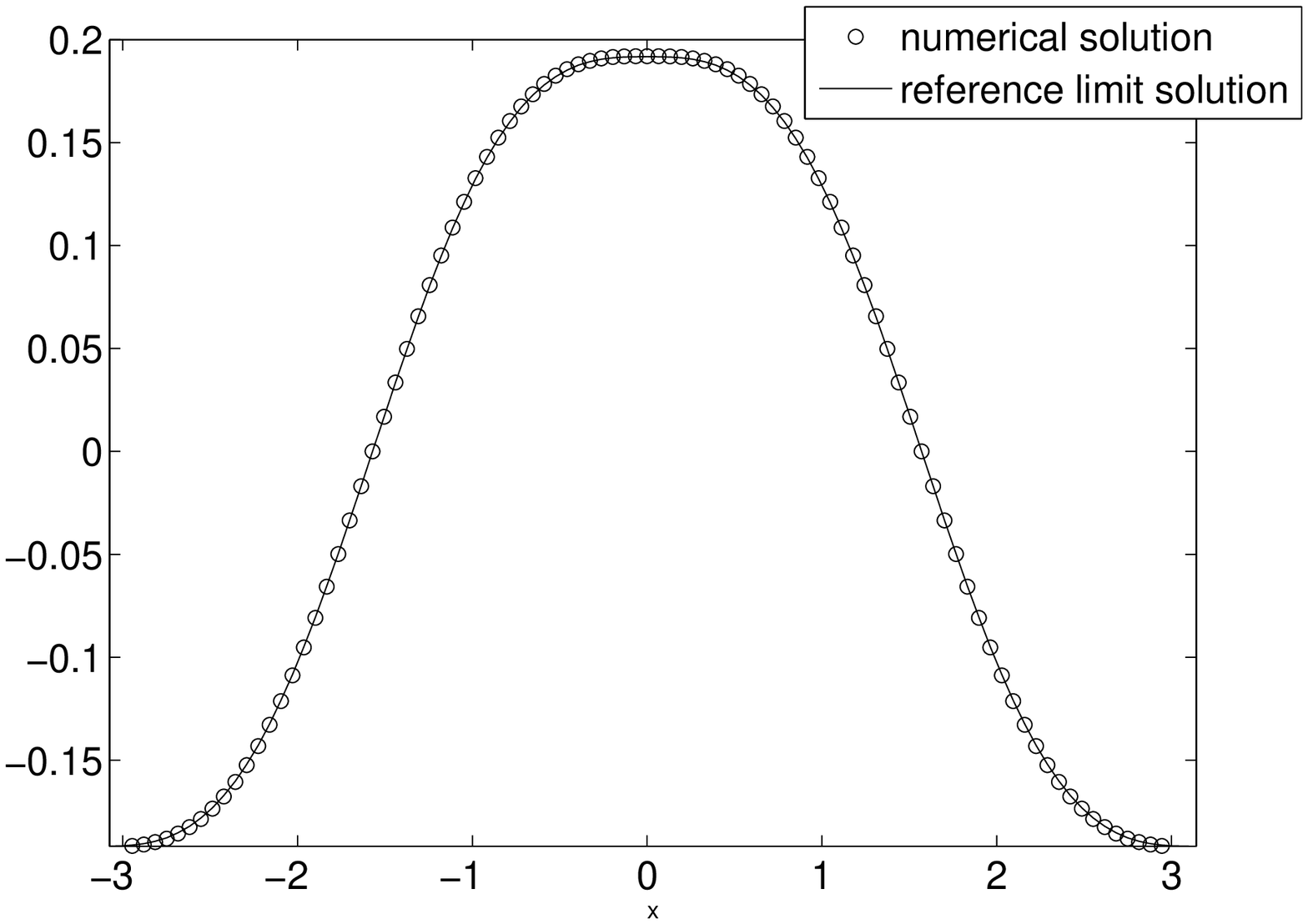}
\caption{ $\circ$ numerical solution with $N = 96$ cells and solid line the reference limit solution with $N = 384$ cells at time $T = 1$. The numerical solutions are compared with $\varepsilon = 10^{-4}$ and $\Delta t = C \Delta x^2$. On the left-hand side, case $m = 0.5$ and $C = 1$.  On the right-hand side, case $m = 2$ and $C = 0.025$.}
\end{figure}

The nonlinear parabolic equation (\ref{eq:redproblem2bis})  has regular solutions if $m \geqslant 1$, i.e. $\alpha \leqslant 0$, while it develops singularities in the derivatives if $0<m<1$, i.e. $\alpha > 0$. Two value of $m$ are considered here, namely $m = 1/2$ ($\alpha = 1$) and $m = 2$ ($\alpha = -1$). The profile of the solution computed with $N = 96$ points is reported in Fig. \ref{Fig:3}. The time step has been chosen to satisfy $\Delta t = C \Delta x^2$, with $C = 0.025$ (right panel with $m = 2$) and $C = 1$ (left panel with $m = 1/2$). The scheme seems to converge in both cases.

Using $N = 12 \cdot 2^p$ with $p = 0,1,...,5$, one obtains the convergence table \ref{tab1}. The $L^{\infty}$ convergence is second order for $m = 2$, and less for $m = 0.5$, as expected from the lack of regularity of the solution. Convergence rate in the latter case improves when measured in slightly weaker norms, i.e. $L^1$ and $L^2$ norm. Observe that second--order is observed even if the scheme is first--order in time, which is due to the parabolic scaling of the time step.

\begin{figure}\label{Fig:3bis}
\centering 
\includegraphics[width=0.85\textwidth]{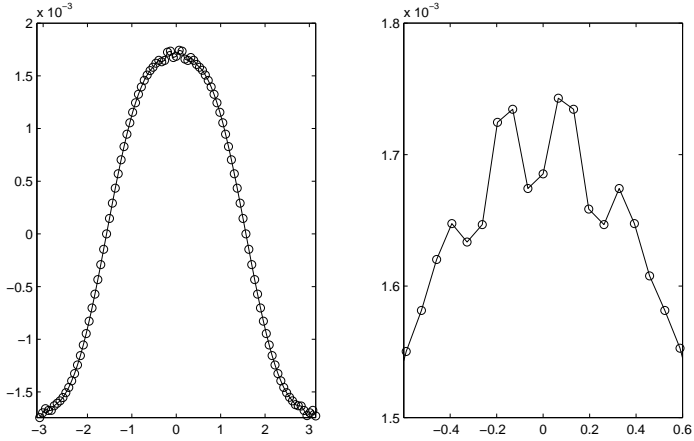}
\caption{On the left-hand side: numerical solution with $N=96$ cells at time $T=1.77$ for $m = 2$, $\varepsilon = 10^{-4}$ and $C = 0.025.$ On the right-hand side: zoom of the region where the instability appears.}
\end{figure}

Let us take a closer look at the case $m = 2$. By integrating for a longer time, some instabilities appear (see Fig. \ref{Fig:3bis}).  The reason of such instabilities is that the following equation (\ref{eq:redproblem2bis}) may be written in the form 
\begin{eqnarray}\label{eq:ceq}
u_t = ((\alpha +1)|u_x|^\alpha) u_{xx}.
\end{eqnarray}
where the low order term $(\alpha +1)|u_x|^\alpha$ plays the role of the diffusive coefficient.

It can be shown \cite{BR} that the first order Implicit-Explicit Euler scheme applied to the linear system (\ref{eq:relax2bis})  suffers from the following stability restriction, independently of $\varepsilon$ 
\begin{eqnarray*}
\frac{\Delta t}{\Delta x^2} \leqslant 1, 
\end{eqnarray*}
which suggests the following condition in the nonlinear case
\begin{eqnarray*}
(\alpha +1)|u_x|^\alpha \, \frac{\Delta t}{\Delta x^2} \leqslant 1. 
\end{eqnarray*}
This condition is used to determine the optimal time step for $m \leqslant 1$, while it shows that no time step can guarantee stability near local extrema if $m  > 1$. For this reason, it is strictly necessary to use {\sl implicit schemes} when $m >1$.

\begin{table}[htbp]
\centering
\subtable[$m = 0.5$]{%
\begin{tabular}{lll}
\multicolumn{1}{c}{\textbf{N}} &
\multicolumn{1}{c}{\textbf{Error}} &
\multicolumn{1}{c}{\textbf{Order}} 
\\
\hline
   12 & 1.2942e-01 & -- \\
   24 & 5.7400e-02  & 1.17\\
   48 & 2.2568e-02 & 1.34\\
   96 & 7.8081e-03 &  1.53\\
   192 & 2.6057e-03&  1.58\\
   384 & 8.2321e-04&  1.66 \\
\hline
\end{tabular}
}\qquad\qquad
\subtable[$ m = 2$]{%
\begin{tabular}{lll}
\multicolumn{1}{c}{\textbf{N}} &
\multicolumn{1}{c}{\textbf{Error}} &
\multicolumn{1}{c}{\textbf{Order}} 
\\
\hline
  12 & 7.9684e-01 & -- \\
   24 & 1.5843e-01  & 2.33\\
   48 & 3.8728e-02 & 2.03\\
   96 &  9.3970e-03& 2.04\\
   192 & 2.3082e-03& 2.02\\
  384 & 5.4599e-04&  2.07\\
\hline
\end{tabular}
}
\caption{$L^{\infty}$-norms of the relative error and convergence rates of $u$ with $C =1$ for case $m = 0.5$ and $C = 0.025$ for $m = 2$.}
\label{tab1}
\end{table}
\begin{table}[htbp]
\centering
\subtable[$L^{2}$-norm]{%
\begin{tabular}{lll}
\multicolumn{1}{c}{\textbf{N}} &
\multicolumn{1}{c}{\textbf{Error}} &
\multicolumn{1}{c}{\textbf{Order}} 
\\
\hline
   12 & 1.8038e-01 & -- \\
   24 & 4.6859e-02  & 1.94\\
   48 & 1.3354e-02 & 1.81\\
   96 & 3.575e-03 &  1.90\\
   192 & 9.3110e-03&  1.94\\
   384 & 2.3020e-04&  2.01 \\
\hline
\end{tabular}
}\qquad\qquad
\subtable[$ L^{1}$-norm]{%
\begin{tabular}{lll}
\multicolumn{1}{c}{\textbf{N}} &
\multicolumn{1}{c}{\textbf{Error}} &
\multicolumn{1}{c}{\textbf{Order}} 
\\
\hline
  12 & 1.9626e-01 & -- \\
   24 & 4.0244e-02  & 2.28\\
   48 & 1.0851e-02 & 1.89\\
   96 &  2.7977e-03& 1.95\\
   192 & 7.0519e-04& 1.98\\
  384 & 1.6909e-04&  2.06\\
\hline
\end{tabular}
}
\caption{$L^{1}$ and $L^{2}$-norms of the relative error and convergence rates of $u$ for $m = 0.5$ with $C = 1$.}
\label{tab2}
\end{table}


\subsection{Removing parabolic stiffness for Kawashima-LeFloch model}

As we have seen, the previous numerical example outlined above converges to an explicit scheme for the limit diffusion equation, and therefore subjected to the classical parabolic CFL restriction $\Delta t \leqslant  C\Delta x^2$. Furthermore in some cases such restriction does not allow the computation of the solution, due to the unboundedness of the diffusion coefficient.

Now, in order to remove the parabolic CFL restriction we adopt a similar technique proposed  for the linear relaxation model that consists in  adding and subtracting the same term to the first equation. In the nonlinear model equation (\ref{eq:relax2}), with $q(u) = 0$, $b(u) = u$, the system becomes
\begin{eqnarray}\label{systempq3} 
  u_t          & = &- (v+\mu(\eps) |u_x|^{\alpha}u_x)_x + \mu(\eps)(|u_x|^{\alpha}u_{x})_x,
   \nonumber        \acapo
               &   &                        \acapo
 \eps^2 v_t & = & -u_x-|v|^{m-1}v. 
        \nonumber
\end{eqnarray}
Here, $\mu(\eps)$ is such that $\mu : \mathbb{R}^{+} \to [0, \ 1]$ and $\mu(0)  =1$. When $\eps$ is not small there is no reason to add and subtract the term $\mu(\eps) u_{xx}$, therefore $\mu(\eps)$ will be small in such a regime, i.e. $\mu(\eps) = 0$. 

We note that in both the BR and BPR approach the term $(|u_x|^{\alpha}u_{x})_x$ is treated implicitly. In the case of a fully implicit treatment of such term, both approaches can be used. Of course the BR approach requires that the Runge-Kutta IMEX is globally stiffly accurate \cite{BPR,BR}.

If, on the other hand, we want to use the treatment of the implicit term by the technique outlined in section  \ref{sec:reflim}, then we are forced to use BPR approach, since the requirement that $b = \tilde{b}$ cannot reconcile with global stiff accuracy.  
To be more specific,  in our case $F(y^{\ast}, y)$ is given by 
$
y = \left(\begin{array}{l}
u\\
v
\end{array}\right),\,
y^{\ast} = \left(\begin{array}{l}
u^{\ast}\\
v^{\ast}
\end{array}\right),
$ and 
$$
F(y^{\ast}, y) = \left(\begin{array}{l}
 -(v^{\ast}_x + \mu(\varepsilon)(|u^{\ast}_x|^{\alpha} u^{\ast}_x)_x) +  \mu(\varepsilon)(|u^{\ast}_x|^{\alpha}u_x)_x \\
-u_x + |v|^{m-1}v 
\end{array}\right).
$$
It can be shown that such semi-implicit R--K scheme is indeed a particular case of IMEX R--K scheme (Pareschi, Russo, private communication). Furthermore, if $\tilde{c} = \tilde{A}e = c = A e$, the classical order conditions on the two schemes 
 \begin{displaymath}
\begin{array}{c|c}
\widetilde{c} & \widetilde{A}\\ 
\hline \\[-.3cm]
 & \ \ \widetilde{b}^T \end{array} \quad
\begin{array}{c|c}
{c} & {A},\\
\hline \\[-.3cm]
 & \ \ {b^T} \end{array}
\end{displaymath}
for order $p = 1,2 ,3$, will guarantee that the semi-implicit R--K has the same order $p$ (cf.~\cite{BFR}).
In order to show the validity of this technique for the computation of the solution for the system (\ref{systempq3}), we consider again the convergence test proposed before. Here we use a classical second order IMEX R--K scheme of type A with $b_i = \tilde{b}_i$, i.e. SSP(3,3,2) {\cite{PaRu}}, for the computation of the numerical solution of the system (\ref{systempq3}). 

As a reference solution we propose again the numerical solution of the limiting equation (\ref{part}) and we compare it with the numerical one computed from the system (\ref{systempq3}). The $L^{\infty}$-convergence is second order for $m = 2$ (cf.~Table \ref{Tab3}) and a CFL hyperbolic condition is used for the time step. 

\begin{figure}
\centering 
\includegraphics[width=0.5\textwidth]{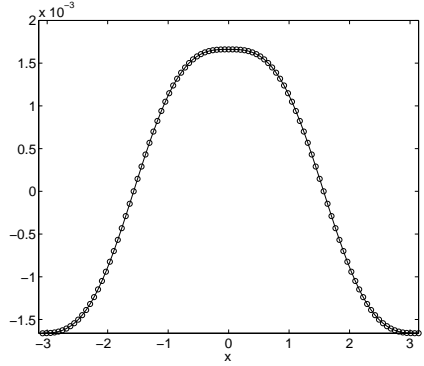}
\caption{Numerical solution with $N=96$ cells at time $T=1.77$ for $m = 2$, $\varepsilon = 10^{-4}$ and $\Delta t = C \Delta x$, with $C = 0.25$.}
\label{Fig:Stab}
\end{figure}

The following remarks are in order: 
\begin{itemize} 
\item In the evaluation of the numerical solution for the case $m = 2$, i.e. $\alpha = -1$, we set a tolerance $TOL$ for computing the derivative $(|u_x|+ TOL)^{\alpha} $ in order to avoid that the derivatives goes to infinity. We chose in our numerical test $TOL = 10^{-12}$.
\item Observe that by integrating for a longer time with a hyperbolic CFL for the case $m = 2$ the numerical solution decreases rapidity to zero without any instabilities (cf.~Fig. \ref{Fig:Stab}). 
\end{itemize}

\begin{table}[t]
\centering
{\begin{tabular}{cccc}
\hline
 &  {\bf N} & {\bf Error}  & {\bf Order}  \\
\hline
 &  12 & 1.6921e-01 &  --\\
                      &  24 & 4.2166e-02  & 2.004\\
                      &  48 &  1.0328e-02& 2.029\\
                      &  96& 2.5371e-03& 2.025\\
                      &192&  6.0394e-04&2.070\\
                      &384 &1.2064e-04 &2.323\\
\hline
\end{tabular}}
\caption{{$L^{\infty}$-norms of the relative error and convergence rates of $u$ with $\Delta t = C \Delta x$ and $C =0.06$ for case $m = 2$.} }
\label{Tab3}
\end{table}
  

\section{Concluding remarks} 

In this paper,
 we have presented a new asymptotic-preserving method in order to construct new IMEX Runge-Kutta schemes which are especially adaoted to deal with 
a class of nonlinear hyperbolic systems containing
 nonlinear diffusive relaxation. As a by-product, effective semi-implicit schemes for the limiting (nonlinear) diffusion equations have been proposed. 
These schemes are able to solve the hyperbolic systems without any restriction on the time step which would classically be imposed by 
the stiff source or by the unboundedness of the characteristic speeds. 
In the limit as the relaxation parameter vanishes, the proposed schemes relax
to implicit schemes for the limit nonlinear convection-diffusion
equation, thus overcoming the classical parabolic CFL condition in the
time step. Although a time discretization up to third-order is available \cite{BPR}, we used here second--order schemes, since  
space discretization is limited to second--order accuracy. The construction of third--order accurate schemes in
space is not a trivial matter in view of the nonlinearity and is by itself an interesting field of investigation. 


\end{document}